\documentclass[a4paper,11pt]{article}
\usepackage{geometry,amssymb,amsmath,enumerate,latexsym,array}
\usepackage[dvips]{graphicx}
\usepackage[all]{xy}
\newcommand{\ep}{\varepsilon} 
\newcommand{\kmods}{\mathsf{KMods}} 
\newcommand{\bA}{\mathsf{A}} 
\newcommand{\bB}{\mathsf{B}}

\newcommand{\sets}{\mathsf{Sets}}

\def\ob{{\mathrm{ Ob}}} 
 
\newcommand{\arr}{\mathrm{Arr}} 
\newcommand{\lan}{\langle} 
\newcommand{\ran}{\rangle} 
\newcommand{\lra}{\leftrightarrow} 

\begin{document}
\newtheorem{example}{Example}[section]
\newtheorem{alg}[example]{Algorithm} 
\newtheorem{rem}[example]{Remark} 
\newtheorem{lem}[example]{Lemma} 
\newtheorem{thm}[example]{Theorem} 
\newtheorem{cor}[example]{Corollary} 
\newtheorem{Def}[example]{Definition} 
\newtheorem{prop}[example]{Proposition} 
\newenvironment{proof}{\noindent {\bf Proof} }{\mbox{} \hfill $\Box$
\mbox{}\\} 
\newenvironment{outline}{\noindent {\bf Outline Proof} }{\mbox{} \hfill $\Box$
\mbox{}\\}

 \title { Gr\"obner Basis Techniques for Computing Actions of K-Categories 
 \thanks{KEYWORDS: Gr\"obner basis, K-category, Action, Kan extension.}}
\author{ Anne Heyworth\thanks{Supported by a University of Wales,
                         Bangor, Research Assistantship.}\\
School of Informatics,\\ University of Wales, Bangor\\ Gwynedd, LL57 1UT, U.K.\\
a.l.heyworth@bangor.ac.uk}

\maketitle

\begin{center} Abstract \end{center}
{\small This paper involves categories and computer science.
Gr\"obner basis theory is a branch of computer algebra which has been
usefully applied to a wide range of problems. Kan extensions are a key concept 
of category theory capable of expressing most algebraic structures. The paper
combines the two, using Gr\"obner basis techniques to compute certain kinds of
Kan extension. 
}

\section{Introduction}

The paper is motivated by a question which arises from two pieces of
research.
Firstly, the work of Brown and Heyworth \cite{paper2}
which extends rewriting techniques to enable the computation of left Kan 
extensions over the category of sets. 
It is well known that left Kan extensions can be defined over categories other 
than $\sets$. 
Secondly, the `folklore', made explicit in \cite{paper3} that rewriting theory 
is a special case of noncommutative Gr\"obner basis theory. 
It is therefore natural to ask whether Gr\"obner bases can provide a 
method for computing Kan extensions beyond the special case of rewriting.\\

To answer this question completely, fully exploiting the computational power
of Gr\"obner basis techniques relating to Kan extensions is the ultimate aim. 
This paper provides a first step by showing how standard noncommutative
Gr\"obner basis procedures can be used to calculate left Kan extensions of
$K$-category actions. In the final section of the paper a number of interesting 
problems arising from the work are identified.

\section{Background}

This paper builds on work of Brown and Heyworth \cite{paper2} on
extensions of rewriting methods.\\

The standard expression of rewriting is in terms of words $w$ in a
free monoid $\Delta^*$ on a set $\Delta$. This may be extended to terms
$x|w$ where $x$ belongs to a set $X$ and the link between $x$ and $w$ is in
terms of an action. More precisely, we suppose a monoid $A$ acts on the
set $X$ on the right, and there is given a morphism of monoids $F$:
$A\to B$ where $B$ is given by a presentation with generating set
$\Delta$.  The result of the rewriting will then be normal forms for
the {\it induced action} of $B$ on $F_*(X)$. This gives an important 
extension of rewrite methods. \\

In fact monoids may be replaced by categories, and sets by directed
graphs. This gives a formulation in terms of left Kan extensions, or
induced actions of categories, which is explained in \cite{paper2}.
Further, categories can be replaced by $K$-categories, as will be
described later.\\

Let $\bA$ be a category. 
  A \textbf{category action} $X$ of $\bA$ is a 
functor $X:\bA \to \sets$. 
Let $\bB$ be a second category and let $F:\bA \to \bB$ be a  functor. 
Then an \textbf{extension of the action $X$ along $F$} is a pair $(E,\ep)$ 
where $E:\bB \to \sets$ is a  functor and $\ep:X \to E \circ F$ is 
a natural transformation.
The \textbf{left Kan extension of the action $X$ along $F$} is an 
extension of the action $(E,\ep)$ with the universal property that 
for any other extension of the action $(E',\ep')$ there exists a 
unique natural transformation $\alpha:E \to E'$ such that 
$\ep'=\alpha \circ \ep$. \\

The problem that has been introduced is that of ``computing a Kan 
extension''. In order to keep the analogy with computation and 
rewriting for presentations of monoids we propose a definition of 
a presentation of a left Kan extension. The papers  
\cite{BLW, CaWa1, CaWa2, Rosebrugh} were very influential on our choices.\\

A \textbf{left Kan extension data} $(X',F')$ consists of small categories
$\bA$, $\bB$ and functors $X':\bA \to \sets$ and $F':\bA \to \bB$.
A \textbf{left Kan extension presentation} is a quintuple
$\mathcal{P}:=kan \lan \Gamma | \Delta | RelB | X | F \ran$ where
$\Gamma$ and $\Delta$ are (directed) graphs; 
$X: \Gamma \to \sets$ and $F: \Gamma \to P \Delta$ are graph 
morphisms to the category of sets and the free category on 
$\Delta$ respectively; 
and $RelB$ is a set of relations on the free category $P\Delta$.

Formally, we say $\mathcal{P}$ \textbf{presents} the left Kan extension 
$(E,\ep)$ of the left Kan extension
data $(X',F')$ where $X':\bA \to \sets$ and $F':\bA \to \bB$ if
 $\Gamma$ is a generating 
graph for $\bA$  and $X: \Gamma \to \sets$ is the restriction of
 $X':\bA \to \sets$;
$cat \lan \Delta | RelB \ran$ is a category presentation for $\bB$ and
$F:\Gamma \to P\Delta$ induces $F':\bA \to \bB$.\\

We expect that a left Kan extension $(E, \ep)$ is given by a set $EB$ 
for each $B \in \ob \Delta$ and a function $Eb: EB_1 \to EB_2$ for 
each $b:B_1 \to B_2 \in \bB$  (defining the functor $E$) together 
with a function $\ep_A: XA \to EFA$ for each $A \in \ob \bA$ (the 
natural transformation). \\

The main result of \cite{paper2} defines rewriting procedures on
$T:= \bigsqcup_{B \in \ob\Delta} \bigsqcup_{A \in 
\ob \Gamma} XA \times P\Delta(FA,B) $
which is basically a set with a partial right action of the arrows of $P \Delta$.\\

Two kinds of rewriting are involved here. The 
first is the familiar $x|ulv \to x|urv$ given by a relation 
$(l,r)$ -- these rules are known as the `$E$-rules'. 
The second derives from a given action of certain words 
on elements, so allowing rewriting $x|F(a) v \to x \cdot a|v$ -- these
rules are known as the `$\ep$-rules'. 
Further, the elements $x$ and $x \cdot a$ may belong to different 
sets. When such rewriting procedures complete, the associated normal form 
gives in effect a computation of what we call the 
{\em Kan extension defined by the presentation}.\\
 
The `folklore' of the relation of rewriting and Gr\"obner basis techniques,
alluded to in  \cite{TMora} and \cite{Birgitsthesis} is made explicit in 
\cite{paper3}.\\

The polynomial ring $K[X^*]$ consists of all polynomials having coefficients
in the field $K$ and terms from $X^*$ together with the usual operations
of polynomial addition and (noncommutative) multiplication.
Given a generating set for an ideal $I$ in this ring it is a problem to
determine whether two given polynomials $f$ and $g$ are equivalent modulo the
ideal ~ i.e. whether they occur within the same congruence class. 
If a Gr\"obner basis $G$ can be constructed for $I$ from the original 
generating set then the congruence problem
can be solved. 
The Gr\"obner basis calculation depends on a well-ordering of $X^*$ and a 
definition of polynomial reduction, which is determined by comparing leading
terms. In the noncommutative case it is not always successful.\\

The key observation is that the rewriting techniques used in calculating a
monoid $M$ from a set of generators $X$ and a rewrite system $R$ compatible
with an ordering $>$ corresponds step-by-step to the Gr\"obner basis 
techniques used in calculating the congruence classes of the polynomial ring 
$K[X^*]$ with respect to the ideal generated by the difference binomials
$l-r$ for $(l,r)$ in $R$.\\

This provides the background to our problem of determining whether
Gr\"obner bases can be used to calculate Kan extensions other than in the
special case of rewriting systems. The first observation is that Gr\"obner
bases involve polynomials, so we should examine how the addition operation
is represented in categories.

\section{$K$-Category Actions}

We use the definitions of \cite{Mitchell}.
Let $K$ be a field.
A \textbf{$K$-category} is a category whose hom-sets (a hom-set is the set of
all
morphisms between a given pair of objects) are $K$-modules.
A morphism of $K$-categories or \textbf{$K$-functor} $F$ preserves the
$K$-module structure of the hom-sets so $F(a+b)=F(a)+F(b), \
F(ka)=kF(a)$
for all arrows $a,b$ such that $a+b$ is defined and scalars $k$ in $K$.\\

The \textbf{free $K$-category on the graph $\Delta$} is the category 
$P_K\Delta$ whose objects are the 
objects of $\Delta$ and whose arrows $\arr P_K\Delta$ are all
polynomials of the form $ p=k_1 m_1 + \cdots + k_n m_n $
where 
$k_1, \ldots, k_n \in K$ and 
$m_1, \ldots, m_n \in P\Delta(B_1,B_2)$ for some $B_1, B_2 \in \ob \Delta$. 
We will refer to $m_1, \ldots, m_n$ as the terms which occur in $f$.
Note that functions $src$ and $tgt$ are well-defined as 
$src(f) := src(m_1) = \cdots = src(m_n)$ and
$tgt(f) := tgt(m_1) = \cdots = tgt(m_n)$.\\

The \textbf{relations of a $K$-category} can be of the form
$p=q$  where both sides have the same source and target.
Therefore $R$ will be assumed to be a set of polynomials $p-q$ ~ 
i.e. a subset of $\arr P_K\Delta$.
If $R=\{r_1,\ldots,r_n\}$ is such a set of relations on $P_K\Delta$ then
the \textbf{congruence generated by $R$} is defined as follows:
$$
f =_R h \text{ if and only if } 
f = h + k_1p_1r_1q_1 + \cdots + k_np_nr_nq_n
$$
for some $k_1, \ldots, k_n \in K$ and 
$p_1, \ldots, p_n, q_1, \ldots, q_n \in \arr P_K\Delta$
where
$src(f) = src(h) = src(u_1) = \cdots = src(u_n)$ and
$tgt(f) = tgt(h) = tgt(v_1) = \cdots = tgt(v_n)$ and
$u_1r_1v_1,\ldots,u_nr_nv_n$ are
defined in $\arr P_K\Delta$.
The $K$-category $P_K\Delta/\! =_R$ whose elements are the congruence
classes of
$\arr P_K\Delta$ with respect to $F$ is known as the \textbf{factor
$K$-category}.

\begin{Def}
Let $K$ be a field. A \textbf{$K$-category presentation} is a pair
$cat_K \lan \Delta | R \ran$ where 
$\Delta$ is a graph and 
$R\subseteq \arr P_K \Delta$.  
The $K$-category it presents is the factor category
$P_K\Delta/\! =_R$.
\end{Def}

Our first result enables the use of Buchberger's algorithm to compute 
Gr\"obner bases which
enable the specification of the morphisms of a $K$-category presented
in this
way.\\

Let $>$ be an admissible well-ordering on $\arr P \Delta$ i.e. $>$ is Noetherian
and compatible with the operation of path concatenation.
Define the \textbf{leading term} of a polynomial $f$ to be the term occurring
in $f$ which is the greatest path in $\Delta$ with respect to $>$ and denote it
$\mathtt{LT}(f)$.
Define a \textbf{reduction relation} $\to_R$ on $\arr P_K\Delta$ by
$f \to f - k_iu_ir_iv_i$ when $u_i(\mathtt{LT}(r_i))v_i$ occurs in $f$ with
coefficient $k_i \in K$ for $u_i,v_i \in \arr P\Delta$, $r_i\in R$.
The reflexive, symmetric and transitive closure of $\to_R$ is
denoted $\stackrel{*}{\lra}_R$.
If the reduction relation $\to_R$ is complete (i.e. Noetherian and
confluent) then we say that $R$ is a \textbf{Gr\"obner basis}.

\begin{lem}
$$
\frac{\arr P_K \Delta}{=_R}\cong\frac{\arr P_K \Delta}{\stackrel{*}{\lra}_R}
$$
\end{lem}

\begin{proof}
It is clear from the definitions that the equivalence relation 
$\stackrel{*}{\lra}_R$ is contained in $=_R$.

For the converse, suppose $f=_R h$. Then there exist
$r_1, \ldots, r_n \in R$ and 
$p_1,\ldots,p_n,q_1,\ldots,q_n\in P_K\Delta$, such that
$f = h + p_1r_1q_1 + \cdots + p_nr_nq_n$. 
By splitting $p_i$ and $q_i$ into their
component terms for $i=1, \ldots, n$ we obtain
$f= h + k_1u_1r_1v_1 + \cdots + k_ju_jr_iv_j + \cdots +  k_tu_tr_nv_t$ for some
$k_1,\ldots,k_t\in K$, $u_1,\ldots,u_t,v_1,\ldots,v_t\in P\Delta$. It
follows immediately from this that $f \stackrel{*}{\lra}_R h$.
\end{proof}

\begin{prop}
The relation $\to_R$ is Noetherian on $\arr P_K\Delta$.
\end{prop}


The \textbf{matches} of $R$ are the pairs of polynomials $(r_1,r_2)$ whose
leading terms overlap on some subword ~ i.e. 
$u \mathtt{LT}(r_1) v = \mathtt{LT}(r_2)$ or 
$\mathtt{LT}(r_2) = u \mathtt{LT}(r_2) v$ or
$u \mathtt{LT}(r_1) = \mathtt{LT}(r_2) v$ or
$\mathtt{LT}(r_1) v = u \mathtt{LT}(r_2)$ for some $u,v \in \arr P\Delta$.
If there is a match between $r_1$ and $r_2$ we may write
$u_1 \mathtt{LT}(r_1) v_1 = u_2 \mathtt{LT}(r_2) v_2$ 
for some $u,v \in \arr P\Delta$.
The \textbf{S-polynomial} resulting from a match is then the difference
$u_1 r_1 v_1 - u_2 r_2 v_2$. The set of S-polynomials of a finite set of 
polynomials is finite and can be computed.

\begin{lem}
If all S-polynomials resulting from matches of $R$ reduce to zero by
$\to_R$ then $\to_R$ is confluent on $\arr P_K\Delta$.
\end{lem}

\begin{outline}
Observing that $\arr P_K\Delta$ is a subset of the free $K$-algebra on
$(\arr P \Delta)^*$ we can deduce that the relation $\to_R$ is confluent on
the free $K$-algebra. 
The fact that  $\to_R$ preserves source and target enables us to deduce that
$\to_R$ cannot reduce an element of $\arr P_K\Delta$ to anything not defined
in $\arr P_K\Delta$. Thus $\to_R$ is confluent on $\arr P_K\Delta$.
\end{outline}


Buchberger's algorithm calculates the S-polynomials of a system $R$ and
attempts to reduce them to zero by $\to_R$. If an S-polynomial cannot be 
reduced it is added to the system. The S-polynomials of the modified system 
$R'$ are then computed -- the process looping until a system is found whose 
S-polynomials can all be reduced to zero.

\begin{thm}[Buchberger's Algorithm and $K$-category Presentations]
If it terminates, then Buchberger's algorithm applied to $(R,>)$, 
will return a Gr\"obner basis for $=_R$ on $\arr P_K \Delta$.
\end{thm}

\begin{proof}
All that remains to be verified is that S-polynomials resulting from
matches found in $R$ can be added to $R$ without altering
$\stackrel{*}{\lra}_R$.
We assume all polynomials in $R$ to be monic (possible since $K$ is a field).
Now S-polynomials result from two types of overlap.

For the first case let $r_1,r_2$ be polynomials in $R$ such that 
$u  \mathtt{LT}(r_1) = \mathtt{LT}(r_2) v$ for some
$u, v \in \arr P \Delta$.
Then the S-polynomial is 
$s := \mathtt{rem}(r_2) v - u \mathtt{rem}(r_1) \in \arr P_K \Delta$ where 
$\mathtt{rem}(r_i) := r_i-\mathtt{LT}(r_i)$ for $i=1,2$. 
Now 
$\mathtt{rem}(r_2) v - u \mathtt{rem}(r_1) = u r_1 - r_2 v$
therefore
$s = \mathtt{rem}(r_2) v - u \mathtt{rem}(r_1) =_R ~0$, 
and hence the congruence generated by $R':=R \cup \{s\}$ coincides with $=_R$.

For the second case let $r_1,r_2$ be polynomials in $R$ such that 
$u \mathtt{LT}(r_1) v = \mathtt{LT}(r_2)$ for some
$u, v \in \arr P \Delta$.
Then the S-polynomial is 
$s := \mathtt{rem}(r_2) - u \mathtt{rem}(r_1) v \in \arr P_K \Delta$. 
Now 
$\mathtt{rem}(r_2) - u \mathtt{rem}(r_1) v = u r_1 v - r_2$
therefore
$s = \mathtt{rem}(r_2) - u \mathtt(r_1) v =_R ~0$, 
and hence the congruence generated by $R':=R \cup \{s\}$ coincides with $=_R$.
\end{proof}

An example of an application of the results proved above can be found
in Section 5.

\section{Left Kan Extensions}

We obtain a further result by expressing the
presentation of a noncommutative polynomial algebra as a problem
of computing a left Kan extension over framed modules $\mathsf{KMods}$
(modules over a fixed field).

\begin{Def}
A \textbf{left Kan extension data for $K$-categories} 
$(M',F')$ consists of small categories
$\bA$, $\bB$ and functors $M':\bA \to \kmods$ and $F':\bA \to \bB$.
A \textbf{left Kan extension presentation for $K$-categories} is a quintuple
$\mathcal{P}:=kan \lan \Gamma | \Delta | RelB | M | F \ran$ where
\begin{enumerate}[i)] 
  \item $\Gamma$ and $\Delta$ are (directed) graphs; 
  \item $M: \Gamma \to \kmods$ and $F: \Gamma \to P_K \Delta$ are graph 
morphisms to the category of $K$-modules and the free $K$-category on 
$\Delta$ respectively; 
  \item  and $RelB$ is a set of relations on the free $K$-category $P_K \Delta$. 
\end{enumerate}
\end{Def}

Formally, we say $\mathcal{P}$ \textbf{presents} the left Kan extension 
$(E,\ep)$ of the left Kan extension
data $(M',F')$ where $M':\bA \to \kmods$ and $F':\bA \to \bB$ if
 $\Gamma$ is a generating 
graph for $\bA$  and $M: \Gamma \to \kmods$ is the restriction of
 $M':\bA \to \kmods$;
$cat_K \lan \Delta | RelB \ran$ is a $K$-category presentation for $\bB$ and
$F:\Gamma \to P_K\Delta$ induces $F':\bA \to \bB$.\\

We expect that a left Kan extension $(E, \ep)$ is given by a set $EB$ 
for each $B \in \ob \Delta$ and a function $Eb: EB_1 \to EB_2$ for 
each $b:B_1 \to B_2 \in \bB$  (defining the $K$-functor $E$) together 
with a function $\ep_A: XA \to EFA$ for each $A \in \ob \bA$ (the 
natural transformation). \\

For the following theorem it is helpful to note that $=^r_{FQ}$ will denote
the right congruence generated by $FQ$. Square brackets $[\cdot]^r_{FQ}$ 
denote the corresponding congruence classes.

\begin{thm}[Congruences on Algebras are Kan Extensions]\mbox{ }\\ 
Let $\mathcal{P}:= kan \lan \Gamma | \Delta | RelB | M | F \ran$
be a presentation of a Kan extension for $K$-categories where:
\begin{enumerate}[i)]
\item $\Gamma$ is the graph with one object $A$ and a collection of arrows $Q$,
\item $\Delta$ is the graph with one object $B$ and a set of arrows $X$, 
\item $RelB$ is a set of polynomial relations $R \subseteq K[X^*]$,
\item $M:\bA \to \mathsf{KMods}$ maps $A$ to $K[1]$ and 
                  the arrows of A to the identity morphism,
\item $F:\bA \to P_K\Delta$ maps the arrows of $\bA$ to polynomials of $K[X^*]$
\end{enumerate}

Then the left Kan extension presented by $\mathcal{P}$ is $(E,\ep)$ where
\begin{enumerate}[i)]
\item $E(B)$ is isomorphic to $(K[X^*]/=_R)/=^r_{FQ}$,
\item $E(b)$ is defined by $E(b)[p]_R:=[pb]_R$,
\item $\ep:M \to E \circ F$ is given by $\ep_AM(q):=[[q]_R]^r_{FQ}$.
\end{enumerate}
\end{thm}

\begin{outline}
It is required to verify that $E$, as defined above, is a well-defined 
$K$-functor. This is quite routine and comes from the fact that the
congruence preserves addition, scalar multiplication and 
right-multiplication. 
To verify that $\ep$ is a natural transformation of $K$-functors
is straightforward, remembering that $M(q)$ is the identity morphism on $K[1]$.
To check the universal property we suppose there is another such pair
$(E',\ep')$ and by drawing the commutative diagram we find that 
there is a unique natural transformation $\alpha:E \to
E'$ defined by
$\alpha(b):=E'(b)(\ep'(1_A))$ for $b \in X^*$.
\end{outline}

It is not claimed that this result is at all deep or difficult, given the
results of \cite{paper2} but it allows the possibility of using 
Gr\"obner bases to compute different types of left Kan extensions.

\begin{cor}
Gr\"obner bases can be used to compute left Kan extensions of the above type.
\end{cor}

\begin{outline}
Let $\mathcal{P}$ be as above. 
Define the $P_K \Delta$-set as
$$
T:= MA \times \arr P_K\Delta
$$ 
and write the terms $A|p$ where $p \in \arr P_K \Delta$.
Define the system of polynomials $\mathcal{S}:=(S_E,S_{\ep})$ where
$$ S_E:= R ~ \text{ and } ~  S_{\ep}:= \{A|Fq-A|1 : q \in Q\}$$

The results in \cite{paper7} describe Gr\"obner
basis procedures for one-sided ideals in finitely presented noncommutative 
algebras over fields. The polynomials defining the $K$-algebra $\bB$ as a 
quotient of the free $K$-algebra $P_K \Delta$ are combined with the 
polynomials defining a right congruence $=^r_{FQ}$ of $\bB$, by using a tagging 
notation. 
Standard noncommutative Gr\"obner basis
techniques can then be applied to the mixed set of polynomials, 
thus calculating $\bB/=^r_{FQ}$ whilst working in a free structure, 
avoiding the complication of computing in $\bB$.\\

Suppose $\mathcal{G}$ is a Gr\"obner basis for $\mathcal{S}$. 
Then the Kan extension is given in the following way:
\begin{enumerate}[i)]
\item $E(B):= \mathtt{IRR}$,
\item $E(b): A|p \mapsto \mathtt{irr}(A|pb)$, for $A|p$ in $E(B)$, $b$ in $K[X^*]$ 
\item $\ep_A(1) := A|1$
\end{enumerate}
where $\mathtt{irr}(A|pb)$ is the irreducible result of repeated reduction of 
$A|pb$ by $\to_{\mathcal{G}}$ and $\mathtt{IRR}$ is the set of all irreducible
terms of $T$.
\end{outline}
 
\begin{rem}
{\em 
It is worth remarking that, as with the rewriting methods developed in 
\cite{paper2}, the Gr\"obner basis methods developed in \cite{paper7}
which are referred to above do not require changes in the existing programs.
The use of tags enables the combination of polynomials giving the conditions
for the action of the Kan extension together with the polynomials giving the 
conditions for the natural transformation.} 
\end{rem}

\section{Examples}

The first example illustrates the previous section, showing that the standard
Gr\"obner basis computation is the computation of a Kan extension and
extending the example to make clear the type of calculation used for
right congruences of algebras.
The second example demonstrates the use of Gr\"obner bases to calculate 
the morphisms of a $K$-category given by a presentation.
In each case we consider the left Kan extension given by a presentation
$\mathcal{P}:=\lan \Gamma | \Delta | RelB | M |F \ran$.

\begin{example}
\emph{Let $\Gamma$ be the trivial graph with one object $A$.
Let $\Delta$ be the graph
with one object $B$ and arrows $X:=[e_1,e_2,e_3]$.
Let $RelB$ be the set of polynomials}
$$
R:=\{e_1e_1-e_1, \ e_2e_2-e_2, \ e_3e_3-e_3, \ e_3e_1-e_1e_3, \
e_2e_1e_2-e_1e_2e_1+\frac{2}{9}e_2-\frac{2}{9}e_1, \
e_3e_2e_3-e_2e_3e_2+\frac{2}{9}e_3-\frac{2}{9}e_2 \}.
$$
\emph{$F:\Gamma \to P_K \Delta$ be inclusion and define $M(A):=K[1]$. 
The system $\mathcal{S}$ consists only of untagged polynomials $R$ because
there are no non-trivial arrows in $\bA$.
We use the length-lexicographic ordering with $e_3>e_2>e_1$ to obtain
Gr\"obner basis for the congruence generated by $\mathcal{S}$ in
$K[X^*]$  by adding} 
$$
e_3e_2e_1e_3-e_2e_3e_2e_1+\frac{2}{9}e_2e_1-\frac{2}{9}e_1e_3
$$
\emph{to $R$.
The irreducible terms $\mathtt{IRR}$ in this case are sums of $K$-multiples of 
the following terms}

\begin{center}
$\{ A|1$,
$A|e_1,\ A|e_2,\ A|e_3$,
$A|e_1e_2,\ A|e_1e_3,\ A|e_2e_1,\ A|e_2e_3,\ A|e_3e_2$,
$A|e_1e_2e_1,\ A|e_1e_2e_3,\ A|e_1e_3e_2,\ A|e_2e_1e_3\ ,A|e_2e_3e_2,\
A|e_3e_2e_1$,
$A|e_1e_2e_1e_3,\ A|e_1e_2e_3e_2,\ A|e_1e_3e_2e_1,\ A|e_2e_1e_3e_2,\
A|e_2e_3e_2e_1\}.$
\end{center}

\emph{In this example the tag ``$A|$'' is redundant: 
the $K$-module $EB$ is a $K$-algebra, in fact it is the Hecke algebra $H_4$. 
Suppose now that $\Gamma$ has one
arrow $q$ whose image under $F$ is $e_2e_1$.
The system of polynomials $\mathcal{S}$ for the Kan extension now has an
$\ep$-polynomial namely $A|e_2e_1 - A|1$. 
Applying Buchberger's Algorithm with the length-lexicographic ordering
$e_3>e_2>e_1>A|$ results in a Gr\"obner basis of mixed polynomials:}

\begin{center}
$\{e_1e_1-e_1, \ e_2e_2-e_2, \ e_3e_3-e_3, \ e_3e_1-e_1e_3, \
e_2e_1e_2-e_1e_2e_1+\frac{2}{9}e_2-\frac{2}{9}e_1, \
e_3e_2e_3-e_2e_3e_2+\frac{2}{9}e_3-\frac{2}{9}e_2 ,$ \\
$e_3e_2e_1e_3-e_2e_3e_2e_1+\frac{2}{9}e_2e_1-\frac{2}{9}e_1e_3, \
A|e_2e_1 - A|1, \ A|e_1 - A|1\}$.
\end{center}

\emph{The right congruence classes of $e_2e_1$ on $H_4$
are represented by sums of $K$-multiples of the following irreducible 
terms ~ i.e. $\mathtt{IRR}$ consists of:} 

\begin{center}
$\{A|1$,
$A|e_2,\ A|e_3$,
$A|e_2e_3,\ A|e_3e_2$,
$A|e_2e_3e_2,\ A|e_3e_2e_1$,
$A|e_2e_3e_2e_1\}.$
\end{center}

\emph{Here the tag ``$A|$'' is necessary in the computation 
of the Gr\"obner basis.
The final results may be written as right congruence classes
$[e_3e_2]^r$, say, instead of tagged terms $A|e_3e_2$ but the representation
as tagged terms allows us to determine whether $e_1e_2e_3$ and $e_2e_3$ occur 
in the same class: reducing $A|e_1e_2e_3$ and $A|e_2e_3$ has the same result, 
so they are congruent.} 
\end{example}

\begin{example}
\emph{Let $\bB$ be the $\mathbb{Q}$-category generated by the graph $\Delta$:}
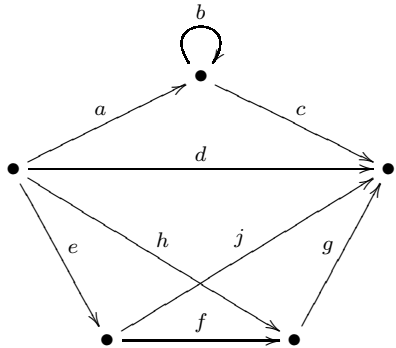
$$
\xymatrix{ && \bullet \ar[drr]^c \ar@(ul,ur)^b &&\\
           \bullet \ar[rrrr]^d \ar[urr]^a \ar[ddr]^e \ar[ddrrr]^h
            &&&& \bullet\\
           &&&&\\
           & \bullet \ar[uurrr]^j \ar[rr]^f && \bullet \ar[uur]^g &\\}
$$
\emph{The arrows of the free category $P_K\Delta$ are sums of 
$\mathbb{Q}$-multiples of terms occuring in the same column of the 
following table (the hom-sets consisting solely of identities are omitted):}
\begin{center}
\begin{tabular}{c|c|c|c|c|c|c|c|c}
$B_1 \to B_2$ & $B_1 \to B_3$ & $B_1 \to B_4$ & 
$B_1 \to B_5$ & $B_2 \to B_2$ & $B_2 \to B_3$ & 
$B_4 \to B_3$ & $B_5 \to B_4$ & $B_5 \to B_3$\\
\hline
$a$ & $ac$ & $h$ &
$e$ & $1_{B_2}$ & $c$ &
$g$ & $f$ & $j$ \\
$ab$ & $abc$ & $ef$ &
 & $b$ & $bc$ &
 & & $fg$ \\
$ab^2$ & $ab^2c$ &&
 & $b^2$ & $b^2c$ &
 & & \\
\vdots & \vdots & &
 & \vdots & \vdots &
 & & \\
$ab^n$ & $ab^nc$ & &
 & $b^n$ & $b^nc$ &
 & & \\  
\end{tabular}
\end{center}
\emph{Let $R$ be the relations defining $\bB$}
$$ R:=\{ab^3-ab^2-ab+a, ~ b^3c-b^2c-bc+c, ~ abc+d-efg, ~ 
        ac+d-hg, ~ fg-j \}
$$
\emph{Applying the length-lexicographic ordering with $a<b<c<d<e<f<g<h<j$
it can be checked that $R$ is a Gr\"obner basis. It can therefore be
immediately deduced that the arrows of $\bB$ are uniquely represented by
$\mathbb{Q}$-multiples of terms occurring in the same column of the following
table:}
\begin{center}
\begin{tabular}{c|c|c|c|c|c|c|c|c}
$B_1 \to B_2$ & $B_1 \to B_3$ & $B_1 \to B_4$ & 
$B_1 \to B_5$ & $B_2 \to B_2$ & $B_2 \to B_3$ & 
$B_4 \to B_3$ & $B_5 \to B_4$ & $B_5 \to B_3$\\
\hline
$a$ & $ac$ & $h$ &
$e$ & $1_{B_2}$ & $c$ &
$g$ & $f$ & $j$ \\
$ab$ & $abc$ & $ef$ &
 & $b$ & $bc$ &
 & & \\
$ab^2$ & $ab^2c$ &&
 & $b^2$ & $b^2c$ &
 & & \\
 & & &
 & \vdots & &
 & & \\
 & & &
 & $b^n$ & &
 & & \\  
\end{tabular}
\end{center}
\end{example}

\section{Further Questions}

\subsection{Induced Modules}
It would be useful to phrase the results of Section 4 in terms of induced 
modules, relating it to the commutative case in \cite{Froberg}.

\subsection{Extensions of Gr\"obner basis techniques}
To apply rewriting to Kan extensions we had to generalise it. We have not yet discovered 
how precisely to generalise Gr\"obner bases to apply to \emph{any} Kan 
extension of $K$-categories over $\mathsf{KMod}$.

\subsection{Rings with Many Objects}
Mitchell's classic work, generalises noncommutative homological ring theory
to (pre)additive category theory \cite{Mitchell}. His work motivates the
investigation of Gr\"obner basis techniques for $K$-categorical Kan 
extensions by the potential for Gr\"obner bases to provide more powerful
methods of computation (of homology or cohomology) in this setting. 

\subsection{Term rewriting and Monads}
Term rewriting systems, widely used throughout computer science, are similar to
algebraic theories (algebraic theories declare term constructors, 
term rewriting
systems declare term constructors and rewrite constructors).
Algebraic theories can be modelled by finitary monads over $\sets$.
Term rewriting systems can be modelled by finitary monads over the
category of preorders $\mathsf{Pre}$.
This has been useful in providing categorical proofs of rewriting theories.
The particularly interesting point is that term rewriting
systems can be modelled as monads over a more complex base category.
So $\mathsf{C}$-algebraic theories can be modelled by finitary monoids on 
$\mathsf{C}$.
There is a relation between monads, adjoint functors and Kan extensions.
We need to investigate the relation between string rewriting for Kan extensions 
and the monads modelling algebraic theories and term rewriting systems.

\subsection{Petri nets}
Gr\"obner basis procedures can be usefully applied in Petri net analysis.
To every Petri net there is an associated category -- a {\em Petri category}
\cite{MandM}.
How does the structure for Petri categories relate to Kan extensions?
Are the Gr\"obner basis techniques usefully extended by relating
these two areas or are they in fact the means by which the areas can be related?

\subsection{Automatic Structure}
For groups, monoids and coset systems there is a well-known concept of an
automatic structure. These systems are special cases of Kan extensions so 
it is natural to ask what would be the definition of an automatic structure 
for a left Kan extension in general.

\end{document}